\documentclass[a4paper,10pt]{article}
\usepackage{amsmath}
\usepackage{amsthm}
\usepackage{amssymb}
\usepackage{latexsym}
\usepackage{enumerate}

\newtheorem{thm}{Theorem}[section]

\title{Infimum of the exponential volume growth and the bottom of the essential 
spectrum of the Laplacian}

\author{Hironori Kumura}

\date{}

\begin{document}

\maketitle

\begin{abstract}
The purpose of this paper is to point out that 
`supremum' in two inequalities of 
Brooks \cite{B1} and \cite{B3} should be replaced with `infimum'. 
\end{abstract}
\section{Introduction}

The Laplace-Beltrami operator $\Delta$ 
on a noncompact complete Riemannian manifold $M$ is 
essentially self-adjoit on $C^{\infty}_0(M)$ 
and its self-adjoit extension to $L^2(M)$ has 
been studied by several authors from various points of view. 
Especially, the bottom $\min \sigma_{{\rm ess}}(-\Delta)$ of 
the essential spectrum of the Laplacian is simply characterized 
by the variational formula
\begin{equation*}
  \min \sigma_{{\rm ess}}(-\Delta)=\lim_{K}\inf_{0\neq f\in C_0^{\infty}(M-K)}
  \frac{\int_M|\nabla f|^2dv_M}{\int_M|f|^2dv_M},
\end{equation*}
where $K$ runs over an increasing set of compact subdomains of $M$ such that 
$\cup K=M$. 
This bottom $\min \sigma_{{\rm ess}}(-\Delta)$ was studied by Donnely, Brooks, and Sunada and so on. 

Donnelly proved in \cite{D} among others that 
$\min \sigma_{{\rm ess}}(-\Delta)\le (n-1)^2k/4$ when 
the Ricci curvature of $M$ is bounded from below by the constant $-(n-1)k$, 
where $n=\dim M$ and $k\le 0$. 
Later, Brooks generalized this Donnelly's 
theorem when the volume of $M$ is infinite:
\begin{thm}[Brooks \cite{B1}]
Let $M$ be a complete Riemannian manifold and set 
$$
    \mu=\limsup_{r\to \infty}
         \frac{ \log {\rm vol}\left(B_{x_0}(r)\right)}{r}.
$$
If the volume of $M$ is infinite, then we have
$$
  \min \sigma_{{\rm ess}}(-\Delta)\le \frac{\mu^2}{4}.
$$
\end{thm}
When the volume of $M$ is finite, Brooks \cite{B3} also proved 
\begin{thm}[Brooks \cite{B3}]
Let $M$ be a complete Riemannian manifold and suppose that the 
volume of $M$ is finite. 
Let us set
$$
   \mu_f=\limsup_{r\to \infty}
   \left[-\frac{\log \left({\rm vol}(M)-{\rm vol}
   (B_{x_0}(r))\right)}{r}\right].
$$
Then we have
$$
   \min \sigma_{{\rm ess}}(-\Delta)\le \frac{\mu_f^2}{4}.
$$
\end{thm}
The purpose of this paper is to improve this two Brooks' theorems.
Theorem $1.1$ is improved as follows, that is, 
`$\limsup$' should be replaced with `$\liminf$':
\begin{thm}
Let $M$ be a complete Riemannian manifold and set 
\begin{align}
    \mu_{{\rm inf}}=\liminf_{r\to \infty}
         \frac{ \log {\rm vol}\left(B_{x_0}(r)\right)}{r}.
\end{align}
If the volume of $M$ is infinite, then we have
$$
  \min \sigma_{{\rm ess}}(-\Delta)\le \frac{\mu_{{\rm inf}}^2}{4}.
$$
\end{thm}
Theorem $1.2$ is also improved as follows, that is, 
`$\limsup$' should be replaced with `$\liminf$':
\begin{thm}
Let $M$ be a complete Riemannian manifold and suppose that the 
volume of $M$ is finite. 
Let us set
$$
   \mu_{f,{\rm inf}}=\liminf_{r\to \infty}
   \left[-\frac{\log \left({\rm vol}(M)-{\rm vol}
   (B_{x_0}(r))\right)}{r}\right].
$$
Then we have
$$
   \min \sigma_{{\rm ess}}(-\Delta)\le \frac{\mu_{f,{\rm inf}}^2}{4}.
$$
\end{thm}
\section{Proof of theorems}
Theorem $1.3$ and $1.4$ will follow the following
\begin{thm}
Let $M$ be a complete Riemannian manifold and $K$ be a compact 
$($possibly empty$)$ set of $M$. 
We denote $\lambda_0(M-K)=\min \sigma(-\Delta _{D,M-K})$, 
where $\Delta _{D,M-K}$ stands for the Dirichlet Laplacian of $M-K$ and 
$\sigma(-\Delta _{D,M-K})$ is the spectrum of $-\Delta _{D,M-K}$. 
We suppose that there exist an increasing sequence $\{K_i\}$ 
of compact subsets of $M$ and positive constants $\alpha$ and $d$ such that 
\begin{align}
  &B_d(K)\subset K_1,\quad \cup_{i=1}^{\infty}K_i=M,\nonumber\\
  &\lim_{i\to \infty}\int_{B_d(\partial K_i)}
  e^{-2\alpha r(x)}dv_M(x)=0.
\end{align}
Here, for $A\subset M$, $B_d(A)$ represents the $d$-neighborhood of $A$ and 
$r(x)$ is the distance function from a fixed point of $M$. 
Then, if 
\begin{align}
  0<\alpha <\sqrt{\lambda_0(M-K)},
\end{align}
we have 
\begin{equation}
   \int_{M-K}e^{2\alpha r(x)}dv_M(x)<\infty.
\end{equation}
\end{thm}
Theorem $2.1$ is proved quite the same way as in Brooks [B1,~Theorem $2$ ], 
and hence, we shall omit its proof. 
\vspace{5mm}

\noindent {\it Proof of Theorem $1.3$}. 
\quad We shall set
\begin{align*}
   &K_t=\{x\in M~|~{\rm dist}(x,K)\le t\}\quad {\rm for}~t\ge 1,\\
   &V(r)={\rm vol}\left(B_r(K)\right).
\end{align*}
Then, the assumption $(1)$ implies that
\begin{align}
  \liminf_{r\to \infty}
         \frac{ \log V(r)}{r}=\mu_{{\rm inf}}.
\end{align}
We remark that if $2\alpha >\mu_{{\rm inf}}$, then there exists an 
increasing sequence $r_i$ of positive numbers such that 
$\lim_{i\to \infty}r_i=\infty$ and 
\begin{align}
     \lim_{i\to \infty}V(r_i+1)e^{-2\alpha r_i}=0.
\end{align}
Indeed, if there is no such sequence $(6)$, 
there exist positive real numbers $\varepsilon$ and $r_0$ such that 
for all $r\ge r_0$ 
$$
     V(r+1)e^{-2\alpha r}\ge \varepsilon.
$$
Hence, we have 
\begin{align}
    \frac{\log V(r)}{r}\ge 2\alpha +\frac{\varepsilon e^{-2\alpha}}{r}
    \quad {\rm for ~all ~}r\ge r_0+1.
\end{align}
But $(7)$ contradicts our assumptions $(5)$ and $2\alpha >\mu_{{\rm inf}}$. 

From $(6)$, we obtain
\begin{align*}
  \int_{B_1(\partial K_{r_i})}e^{-2\alpha r(x)}dv_M(x)
 &\le \left(V(r_i+1)-V(r_i-1)\right)e^{-2\alpha r_i}\\
 &\le V(r_i+1)e^{-2\alpha r_i}\to 0\quad (i\to \infty).
\end{align*}
Thus, when we set $K_i=K_{r_i}$ in Theorem $2.1$, the assumption $(2)$ holds. 
Therefore, we now conclude from Theorem $2.1$ that 
$2\alpha >\mu_{{\rm inf}}$ and $(3)$ imply that $(4)$. 
But it is impossible, since $M-K$ has infinite volume. 
Hence there is no such $\alpha$, and we have 
$2\sqrt{\lambda_0(M-K)}\le \mu_{{\rm inf}}$, that is, 
$\lambda_0(M-K)\le \mu_{{\rm inf}}^2/4$. 
Taking the limit over arbitrary large $K$, we get 
$\min \sigma_{{\rm ess}}(-\Delta)\le \mu_{{\rm inf}}^2/4$. 
We have thus proved Theorem $1.3$.

\vspace{5mm}

\noindent {\it Proof of Theorem $1.4$}. 
\quad For simplicity, we shall set 
$V(r):={\rm vol}(M)-{\rm vol}\left(B_{x_0}(r)\right)$, 
and take a compact subset $K$ of $M$. 
Then, $\lim_{r\to \infty}V(r)=0$, since $M$ has finite volume. 
We remark that if $2\alpha>\mu_{f,{\rm inf}}$, 
then there exists a sequence of positive numbers $\{r_i\}_{i=1}^{\infty}$ 
such that 
\begin{align}
   &B_{x_0}(r_1)\supset K,\nonumber \\
   &r_{i+1}\ge r_{i}+1,\\
   &V(r_{i+1})e^{2\alpha r_i}\le 1,\\
   &2\alpha > \frac{-1}{r_i}\log V(r_i).
\end{align}
The inequality $(9)$ comes from the fact that $\lim_{r\to \infty}V(r)=0$, and 
$(10)$ implies 
\begin{align}
   e^{2\alpha r_i}V(r_i)\ge 1 \quad {\rm for~ all~~}i\ge 1.
\end{align}
Therefore, for any integer $k\ge 1$, we have
\begin{align*}
       &\int_{M-K}e^{2\alpha r(x)}dv_M(x)\\
   \ge & \sum_{i=1}^{k}
         \left\{{\rm vol}(B_{x_0}(r_{i+1}))-{\rm vol}(B_{x_0}(r_i))\right\}
         e^{2\alpha r_i}\\
   =   & \sum_{i=1}^{k}
         \left\{V(r_i)-V(r_{i+1})\right\}
         e^{2\alpha r_i}\\
   =   &V(r_1)e^{2\alpha r_1}
        +\sum_{i=2}^{k}V(r_i)e^{2\alpha r_i}
        \{1-e^{2\alpha(r_{i-1}-r_i)}\}-V(r_{k+1})e^{2\alpha r_k}\\
   \ge   &V(r_1)e^{2\alpha r_1}
        +(k-1)\{1-e^{-2\alpha}\}-1.
\end{align*}
In the last line, we have used $(11)$, $(8)$, and $(9)$ in turn. 
Letting $k\to \infty$, we get
\begin{align}
    \int_{M-K}e^{2\alpha r(x)}dv_M(x)=\infty
\end{align}
Now, when we set $K_i=\{x\in M~|~{\rm dist}(x,K)\le r_i\}$, $K_i$ and $\alpha$ 
satisfy
$$
    \cup_{i=1}^{\infty}K_i=M,~~{\rm and}~~
  \lim_{i\to \infty}\int_{B_1(\partial K_i)}
  e^{-2\alpha r(x)}dv_M(x)=0,
$$
since the volume of $M$ is finite. 
Hence, if $\alpha>\mu_{f,{\rm inf}}/2$ satisfies 
$0<\alpha <\sqrt{\lambda_0(M-K)}$, Theorem $2.1$ implies 
$$
   \int_{M-K}e^{2\alpha r(x)}dv_M(x)<\infty.
$$
But this contradicts $(12)$. 
Therefore, there is no such $\alpha$, and hence, 
we get $\mu_{f,{\rm inf}}^2/4\ge \lambda_0(M-K)$. 
Taking the limit over arbitrary large $K$, we get 
$\min \sigma_{{\rm ess}}(-\Delta)\le \mu_{f,{\rm inf}}^2/4$. 
We have thus proved Theorem $1.4$. 
\section{Example and remark}

In this section, we shall consider the sharpness of our theorems. 
We will begin our discussion by 
considering a rotationally symmetric manifold 
$({\bf R}^n,g=dr^2+f(r)^2g_{S^{n-1}(1)})$. 
Here, we take an increase sequence of positive numbers 
$0<a_1<a_2<\cdots \to \infty$ and define
\begin{align*}
f(r)=
\begin{cases}
     1,~\quad \mathrm{if}~~a_{4k+1}\le r\le a_{4k+2},\\
     e^r,\quad \mathrm{if}~~a_{4k+3}\le r\le a_{4k+4},
  \end{cases}
\end{align*}
where $k=0,1,2,\cdots $. 
We also assume that $f$ is monotone on the intervals $[a_{4k+2},a_{4k+3}]$ and 
$[a_{4k+4},a_{4k+5}]$ for $k=0,1,2,\cdots $. 
Then, if we choose `exponential-intervals' 
$[a_{3}, a_{4}],~[a_{7}, a_{8}],\cdots$, successively large, then 
the function $r^{-1}\log {\rm vol}(B_0(r))$ will oscillate between the values 
$0$ and $n-1$, where $B_0(r)$ is the ball centered at the origin $0$ 
with radius $r$. 
Thus we will have 
\begin{align*}
  &\mu=
  \limsup_{r\to \infty}\frac{\log {\rm vol}\left(B_{0}(r)\right)}{r}=n-1,\\
  &\mu_{{\rm inf}}
  =\liminf_{r\to \infty}\frac{\log {\rm vol}\left(B_{0}(r)\right)}{r}=0,\\
  &\min \sigma_{{\rm ess}}(-\Delta)=0.
\end{align*}
This example shows that Theorem $1.3$ is indeed sharper than Theorem $1.1$. 
We can also construct an example with similar nature which shows that 
Theorem $1.4$ is indeed sharper than Theorem $1.2$. 

On the other hand, Theorem $1.3$ may fail to be sharp. 
Indeed, as is pointed out by Brooks \cite{B1}, 
there exists a solvable group $G$ with exponential growth. 
One such group is given in Milnor \cite{M}. 
Then if $N$ be any compact manifold with $\pi_1(N,x_0)=G$ and 
$M$ is the Riemannian universal cover of $N$, 
a lemma of Milnor says that $\mu_{{\rm inf}}>0$, 
while a Brooks' theorem in \cite{B2} (see also Sunada \cite{S}) implies that 
$\min \sigma_{{\rm ess}}(-\Delta)=0$, since $G$ is amenable. 
The reason why this gap occurs is that the distance spheres need not to be 
the most efficient candidates for the isoperimetric inequalities of $M$. 
Now let us recall the following F$\o$lner-Brooks theorem:
\begin{thm}[F$\o$lner-Brooks \cite{B2}]
Let $N$ be an $n$-dimensional compact Riemannian manifold, 
and $M$ its Riemannian universal cover of $N$. 
Then $\pi_1(N,x_0)$ is amenable if and only if, 
for every $($possibly disconnected$)$ fundamental set of $M$, and for every 
$\varepsilon>0$, there exists a finite subset $E$ of $\pi_1(N,x_0)$ such that 
$$
    H=\bigcup_{g\in E}g\cdot F
$$
satisfies the isoperimetric inequality:
\begin{align*}
   \frac{{\rm vol}_{n-1}(\partial H)}{{\rm vol}_{n}(H)}<\varepsilon.
\end{align*}
\end{thm}
From the proof of Theorem $3.1$ and a lemma of Milnor \cite{M}, we see that 
for any $r_0\ge {\rm diam}\,(N)$, there exist a sequence of 
finite subsets $E_i$ of $\pi_1(N,x_0)$ such that 
$$
    H_i=\bigcup_{g\in E_i}g\cdot B_{x_0}(r_0)
$$
satisfies
\begin{align*}
   \lim_{i\to \infty}\frac{{\rm vol}_{n-1}(\partial H_i)}{{\rm vol}_{n}(H_i)}
   =0.
\end{align*}
Thus, in this case, we see that the distance spheres in Theorem $1.3$ should be replaced with the {\it finite union $H_i$ of distance spheres transformed by 
covering transformations $E_i$}. 

Our main concern above has been the bottom of the essential spectrum. 
But as for the essential spectrum itself, we note that there is a simple 
criterion:
\begin{thm}[\cite{K1},\cite{K2}]
    Let $M$ be a complete Riemannian manifold and assume that there exists an 
    open subset $U$ of $M$ with $C^{\infty}$ compact boundary $\partial U$ 
    such that the outward exponential map 
    $\exp^{\perp}_{\partial U}:N^{+}(\partial U)\to M-U$ 
    induces a diffeomorphism. 
    We set $r(x)={\rm dist}(x,U)$ for $x\in M-U$. 
    If $\Delta r\to c$ as $r\to \infty$ for a constant 
    $c\in {\bf R}$, then $[c^2/4,\infty)\subset \sigma_{{\rm ess}}(-\Delta)$. 
    Moreover, when $U$ is relatively compact, the equality holds$:$
    $\sigma_{{\rm ess}}(-\Delta)=[c^2/4,\infty)$. 
\end{thm}

\vspace{1mm}
\begin{flushleft}
Hironori Kumura\\ 
Department of Mathematics\\ 
Shizuoka University\\ 
Ohya, Shizuoka 422-8529\\ 
Japan\\
E-mail address: smhkumu@ipc.shizuoka.ac.jp
\end{flushleft}

\vspace{20mm}

After submitting this article to this preprint server, Professor Yusuke Higuchi let me know that Theorem $1.3$ and $1.4$ are already known by him (see the following paper). So, I want to delete this article. However, the organizers of this preprint server refuse to delete this article. Hence, I should announce that this article will be unpublished. 

Yu.HIGUCHI,  A remark on exponential growth and the spectrum of the Laplacian, Kodai Math. J. 24 (2001), 42--47.

\end{document}